# Anticipating critical transitions in multidimensional systems driven by time- and state-dependent noise

Andreas Morr,[1,2,*] Keno Riechers,[3] Leonardo Rydin Gorjão,[4] and Niklas Boers[1,2,5]
[1]*Earth System Modelling, School of Engineering and Design, Technical University of Munich, 85521 Ottobrunn, Germany*
[2]*Research Domain IV – Complexity Science, Potsdam Institute for Climate Impact Research, 14473 Potsdam, Germany*
[3]*Climate Physics, Max Planck Institute for Meteorology, 20146 Hamburg, Germany*
[4]*Faculty of Science and Technology, Norwegian University of Life Sciences, 1433 Ås, Norway*
[5]*Global Systems Institute and Department for Mathematics, University of Exeter, EX4 4QF Exeter, United Kingdom*

Anticipating bifurcation-induced transitions in dynamical systems has gained relevance in various fields of the natural, social, and economic sciences. Before the annihilation of a system's equilibrium point by means of a bifurcation, the system's internal feedbacks that stabilize the initial state weaken and eventually vanish, a process referred to as critical slowing down (CSD). In one-dimensional systems, this motivates the use of variance and lag-1 autocorrelation as indicators of CSD. However, the applicability of variance is limited to time- and state-independent driving noise, strongly constraining the generality of this CSD indicator. In multidimensional systems, the use of these indicators is often preceded by a dimension reduction in order to obtain a one-dimensional time series. Many common techniques for such an extraction of a one-dimensional time series generally incur the risk of missing CSD in practice. Here, we propose a data-driven approach based on estimating a multidimensional Langevin equation to detect local stability changes and anticipate bifurcation-induced transitions in systems with generally time- and state-dependent noise. Our approach substantially generalizes the conditions under which CSD can reliably be detected, as demonstrated in a suite of examples. In contrast to existing approaches, changes in deterministic dynamics can be clearly discriminated from changes in the driving noise using our method. This substantially reduces the risk of false or missed alarms of conventional CSD indicators in settings with time-dependent or multiplicative noise. In multidimensional systems, our method can greatly advance the understanding of the coupling between system components and can avoid risks of missing CSD due to dimension reduction, which existing approaches suffer from.

## I. INTRODUCTION

A mechanistic understanding of complex high-dimensional physical systems is essential for assessing the risk of abrupt regime shifts, for example, in ecological, climatic, social, or financial systems. Such shifts may occur when critical forcing thresholds, which correspond to underlying bifurcation points, are crossed [1–4]. Reducing complex systems to a low-dimensional summary observable $\mathbf{X}_t$ has leveraged impressive modeling capabilities [5–7]. This is particularly important because observations are typically available in the form of multivariate time series of just a few dimensions.

Commonly, the dynamics of the summary observable $\mathbf{X}_t \in \mathbb{R}^n$ are approximately separated into a deterministic component $A(\mathbf{X}_t, t)dt$ and a stochastic component $B(\mathbf{X}_t, t)d\mathbf{W}_t$ that represents the action of the omitted dimensions. This results in the Langevin equation

$$d\mathbf{X}_t = A(\mathbf{X}_t, t)dt + B(\mathbf{X}_t, t)d\mathbf{W}_t, \quad (1)$$

commonly employed as a model for partially resolved dynamics. Even though the stochastic component can take more complicated forms, we restrict ourselves to Gaussian white noise $d\mathbf{W}$ and refer to existing extensions to the case of correlated noise [8,9].

This framework facilitates a mathematical description of abrupt regime shifts in terms of dynamic bifurcations in low-dimensional dynamical systems [2,10]. Prior to the transition, the deterministic drift $A(\mathbf{X}_t, t)$ exhibits negative feedback mechanisms keeping the system close to a stable equilibrium $\mathbf{x}_0^*(t)$ with $A(\mathbf{x}_0^*(t), t) = 0$ [11–13]. The explicit time dependence of $A(\mathbf{X}_t, t)$ reflects changes in the forcings that act on the system from the outside and alter the deterministic, coarse-grained dynamics. At the bifurcation point, i.e., at the critical level of forcing, the currently occupied equilibrium state is annihilated and the system abruptly transitions to another stable state $\mathbf{x}_1^*(t)$.

It is well known that a weakening of the negative feedback precedes the annihilation of an attracting equilibrium point [2,14]. Heuristically speaking, this results in weaker and slower responses to the pseudorandom perturbations stemming from the unresolved dynamics. This phenomenon is referred to as critical slowing down (CSD), and it can manifest in an increase of variance and lag-1 autocorrelation (AC(1)) in one or more components of the observable [2,15–17]. These two quantities are therefore often employed to anticipate bifurcation-induced abrupt transitions, and their simultaneous increase has been suggested as an early warning signal (EWS) [2,3,18].

---

*Contact author: andreas.morr@tum.de

We first introduce the notion of local system stability in a one-dimensional dynamical system denoted by

$$dX_t = a(X_t, t)dt + b(X_t, t)dW_t. \quad (2)$$

Linearizing the deterministic equilibrium dynamics yields the linear restoring rate

$$\lambda(t) := -\partial_x a(x^*(t), t). \quad (3)$$

In the annihilation of an equilibrium point, $\lambda$ vanishes, motivating its analysis in the context of anticipating critical transitions. It can be shown that variance and AC(1) of the observable in this linear framework with time- and state-independent noise $b(x, t) \equiv \sigma$ are given by

$$\text{Var}[X] = \frac{\sigma^2}{2\lambda} \xrightarrow{\lambda \to 0} \infty,$$

$$\text{AC}_X(1) = \exp(-\lambda \Delta t) \xrightarrow{\lambda \to 0} 1, \quad (4)$$

where $\Delta t > 0$ is the sampling time step of the data. The extension of the concept of local stability to the general, higher-dimensional setting of Eq. (1) is discussed in the Methods section.

Detection of CSD via variance and AC(1) is usually preceded by a reduction of the system to a one-dimensional observable, either by leveraging physical understanding or employing principal component analysis (PCA) in order to identify a linear combination of components which may be experiencing stability loss [15,19]. In such a reduction, crucial information about system stability may be lost [20,21], and thus methods that avoid deleterious preprocessing are advantageous.

The derivation of the variance and AC(1) in Eq. (4) hinges on the *a priori* assumption of entirely linear feedback. In applications, the system might explore parts of the state space where nonlinearities in the feedback are not negligible anymore, placing the usefulness of Eq. (4) into question. In contrast, the linear restoring rate $\lambda$ directly captures the desired information of local stability at the equilibrium point and should therefore be considered the key quantity to measure system stability and detect CSD. This approach has been adopted for one-dimensional systems, e.g., in [9,22–24].

Time- or state-dependent driving noise can lead to both false negative and false positive EWS [8,25,26]. In particular, variance is vulnerable to alterations in noise levels. In contrast, AC(1) reliably increases during CSD, even under the influence of changing noise strength, provided that the noise is uncorrelated. Given that in real-world applications, the assumption of time- and state-independent noise is hardly justifiable, understanding the evolution of the diffusion term $b(X_t, t)dW_t$ is crucial for assessing stability changes derived from data. In particular, a methodology is needed to extract from an observable a more holistic picture of both the deterministic dynamics and the driving noise.

To obtain an estimation $\widehat{\lambda}$, we first estimate the drift function $a(x)$. In a second step, we perform a spatially local linear fit to $\widehat{a}(x)$ [27,28] in a selected neighborhood around $\widehat{x}^* = \widehat{\mathbb{E}}[X]$ (see Methods and Appendix A). We include an estimation of the diffusion $b(x, t)$ to supplement the standard CSD indicators and to avoid false positives and false negatives caused by changes in the driving noise. In particular, we discuss situations where the conventional CSD indicators give an ambiguous or misleading picture. We show how the method proposed herein conclusively resolves these ambiguities.

We carry the concept to multiple dimensions to obtain a more generally applicable indicator of CSD that is robust with respect to time- and state-dependent noise. Conventional approaches employing PCA to infer information about system stability may fail in such settings.

## II. METHODS

The drift and diffusion coefficients $A$ and $B$ of the $n$-dimensional stochastic differential equation (SDE) (1) have the following representation in terms of the increments $\Delta \mathbf{X}_t := \mathbf{X}_{t+\Delta t} - \mathbf{X}_t$ of the process $\mathbf{X}$ [27,29–35]:

$$A(\mathbf{x}, t) = \lim_{\Delta t \to 0} \frac{1}{\Delta t} \mathbb{E}[\Delta \mathbf{X}_t | \mathbf{X}_t = \mathbf{x}], \quad (5)$$

$$BB^\top(\mathbf{x}, t) = \lim_{\Delta t \to 0} \frac{1}{\Delta t} \mathbb{E}[\Delta \mathbf{X}_t \Delta \mathbf{X}_t^\top | \mathbf{X}_t = \mathbf{x}]. \quad (6)$$

The former relation generally holds for Markovian processes with a drift term $A$, while the latter relation specifically concerns the Gaussian white noise model.

If the SDE in Eq. (1) exhibits effective time independence, i.e., $A(\mathbf{x}, t) \equiv A(\mathbf{x})$ and $B(\mathbf{x}, t) \equiv B(\mathbf{x})$ in some observation time span, we can estimate the stationary variance and AC(1) of each component as well as $A(\mathbf{x})$ and $BB^\top(\mathbf{x})$ through a law of large numbers estimator. More concretely, if the sample path of $\mathbf{X}$ is available at sufficiently small time steps $\Delta t > 0$, one may estimate $A(\mathbf{x})$ and $BB^\top(\mathbf{x})$ by replacing the above ensemble average by the mean of the observed increments. This yields consistent estimators that converge to the true $A$ and $BB^\top$, omitting here a known bias that stems from the nonzero $\Delta t$ [32]. We will see that this bias can incur an unfavorable uncertainty in the estimator $\widehat{\lambda}$.

We generalize the definition of local system stability $\lambda$ in Eq. (3) from the one-dimensional setting to the multidimensional setting of Eq. (1). Consider the Jacobian matrix of $A$ at an equilibrium point $\mathbf{x}^* \in \mathbb{R}^n$ with $A(\mathbf{x}^*, t) = \mathbf{0}$. Such an equilibrium point is stable if and only if all the eigenvalues' real parts $(-\lambda_k + i\omega_k)_{k=1,\ldots,n}$ of the Jacobian matrix $DA(\mathbf{x}^*, t)$ embody negative feedback $-\lambda_k < 0$. Accordingly, we regard the set of $(\lambda_k)_{k=1,\ldots,n}$ as a measure of system stability. When an equilibrium point vanishes by means of a bifurcation, at least one of the $\lambda_k$ passes through zero.

We obtain an estimation of $A$ and $BB^\top$ from time series data using the respective estimators given in [28,36], where kernel-based methods are employed to leverage the law of large numbers convergence. Next, a multivariate ordinary least-squares regression is performed to extract an estimate of the matrix $DA$ (see Appendix A for details). The regression is restricted to a neighborhood of the estimated equilibrium point in order to disregard data points that potentially lie outside the linear regime. On the one hand, this should give a more accurate estimate of the local linear stability compared to conventional methods. On the other hand, such a restriction effectively reduces the number of data points contributing to the estimation and can thus incur larger uncertainties. The real parts of the eigenvalues derived from the estimation of $DA$ then serve as a measure of stability.

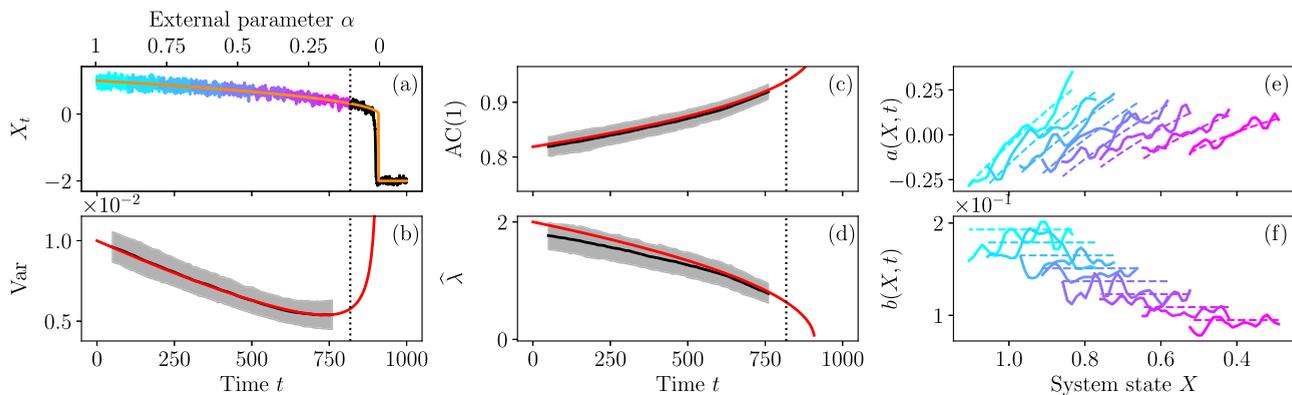

FIG. 1. Application of CSD indicators for synthetic data generated by the model in Eq. (7). The noise strength $\sigma$ is ramped linearly from 0.2 to 0.06 over the integration time span. (a) Sample paths for zero noise (orange) and with noise (divided into colored time series windows). (b), (c) Conventional CSD indicators of variance and AC(1) (black) are calculated on detrended windows and plotted along with the theoretical values (red) obtained from the time-local Ornstein–Uhlenbeck linearization. The black lines constitute ensemble means, and the shaded bands represent the 68% confidence intervals on $N = 1000$ samples. (d) Estimator $\widehat{\lambda}$ as proposed in this work, along with the true value $\lambda = -\partial_x a(x^*(t), t)$ (red). (e), (f) Estimated functions $a(x, t) = -x^2 + \alpha(t)$ and $b(x, t) = \sigma(t)$ plotted in their $x$ dependence. The time dependence is illustrated by the respective color of the plot corresponding to the investigated time series window from (a). The dashed lines represent the true functions $a$ and $b$ at the respective window-center time. The estimated slope of $a$ is taken to be $-\widehat{\lambda}$, which can be seen to diminish over the observation time span. For more details on these estimation routines, see the Methods section or Appendix A. The apparent time dependence of $b$ helps reconcile the inconsistent signals in variance and AC(1) with respect to CSD that may otherwise lead to a missed alarm. All estimations were performed on running windows of length $10^2$ time units consisting of $10^3$ data points, considering the sample time step is $\Delta t = 10^{-1}$.

To investigate the evolution of stability in an application setting, one can assess the eigenvalues of $DA$ locally in time by partitioning a recorded time series into windows. The window length should be short enough to ensure approximately time-independent dynamics within individual windows and long enough to include sufficient data points for robust estimations of drift and diffusion. Negative trends in any of the $\lambda_k$ computed within these sliding windows would indicate a destabilization of the equilibrium state, which may point to an upcoming abrupt transition.

We show that the estimation of the diffusion matrix $BB^\top$ may in certain situations explain peculiar behavior in the conventional CSD indicators, such as simultaneously decreasing variance and increasing AC(1). If instead of local stability, here measured by $\widehat{\lambda}$, one wishes to examine mean exit times from stable equilibria [37], the estimated diffusion matrix $BB^\top$ provides important additional information [38].

## III. RESULTS

### A. Normal-form fold bifurcation with time-dependent noise

The prototypical structure employed for conceptualizing abrupt transitions in many natural systems is the fold bifurcation [2,39]. Consider therefore the SDE defined by

$$dX_t = \left(-X_t^2 + \alpha(t)\right)dt + \sigma(t)dW_t, \quad (7)$$

where $\alpha$ is the bifurcation parameter and $\sigma$ is the noise strength. For positive $\alpha > 0$, there exists a stable equilibrium at $x^*(t) = \sqrt{\alpha}$ which vanishes at the critical threshold $\alpha_{\text{crit}} = 0$ [see Fig. 1(a)]. While decreasing the parameter $\alpha$ linearly, we simultaneously ramp down the noise strength $\sigma(t)$ from 0.2 to 0.06. Such an evolution should be understood as a change in the nature of the omitted fast dynamics, which cannot be ruled out in many applications [40,41].

The temporal evolution of the variance and AC(1) can be approximated after linearizing the system around the time-dependent equilibrium $x^*(t)$ [red lines in Figs. 1(b) and 1(c)]: $\text{Var}[X_t] \approx \sigma(t)^2/4\sqrt{\alpha(t)}$ and $\text{AC}_X(1) \approx \exp(-2\sqrt{\alpha(t)}\Delta t)$. In some cases, the time-dependent noise thus induces a deceiving downward trend in the estimated variance of the system [Fig. 1(b)] alongside increasing AC(1) [Fig. 1(c)]. The conflicting indications given by variance and AC(1) would mislead the observer to conclude that no significant EWS is present. The estimation of the linearized feedback $\lambda$ [Fig. 1(d)] clearly indicates a weakening of local system stability and thus the presence of CSD.

The apparently inconsistent results of the conventional CSD indicators can be understood and reconciled by examining the structure of the drift and diffusion coefficients $a(x, t) = x^2 - \alpha(t)$ and $b(x, t) = \sigma(t)$ and their evolution in time. The true quantities for the functions $a(x, t)$ and $b(x, t)$ at different times $t$ are plotted in Figs. 1(e) and 1(f), along with the estimations obtained during the procedure outlined in the methods section. After disclosing the time dependency of the diffusion coefficient in Fig. 1(f), the diverging trends in variance and AC(1) can be correctly interpreted in a CSD assessment. The decrease in variance can be attributed to the decreasing noise strength, and the approach of a bifurcation can be confirmed. In an analogous setting featuring no bifurcation but an increasing noise strength, the incurred increase in variance can be attributed correctly, and concerns that the system is approaching a bifurcation can be unambiguously dispelled.

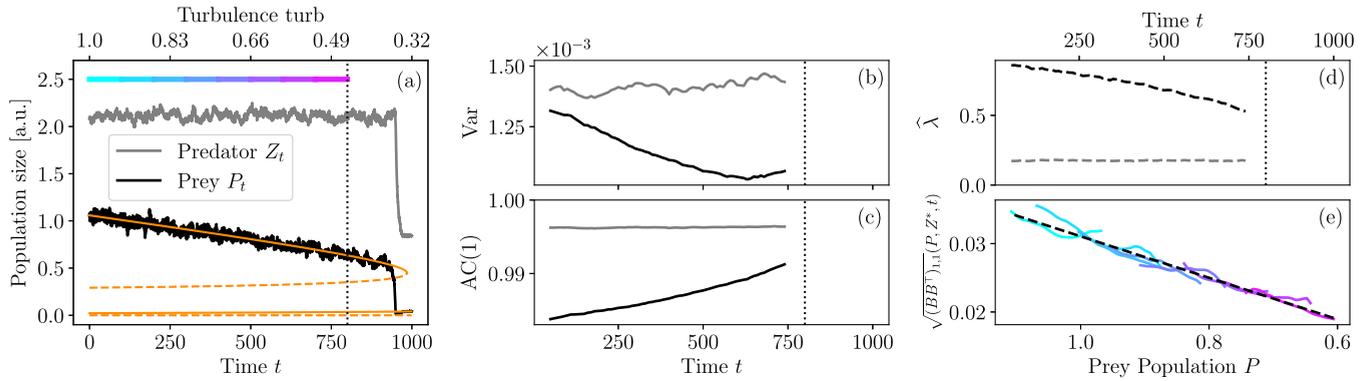

FIG. 2. Application of the CSD indicators on time-series data obtained from the predator-prey model defined in Eqs. (8) and (9). (a) Sample paths of the prey population $P$ (black) and the predator population $Z$ (gray). The stable and unstable equilibria of $P$ in their dependence on turb($t$) are plotted in orange. (b), (c) Means of the conventional CSD indicators variance and AC(1) over $N = 1000$ sample trajectories. (d) Real parts $-\lambda_{1,2}$ of the estimated eigenvalues of the local Jacobian matrix. The eigenvalues have been assigned (by color) to the two populations, as the corresponding eigenspace basis aligns very well. (e) Square root of the estimated matrix element $(BB^\top)_{1,1}(P, Z, t)$ of the single sample trajectory shown in (a). It is plotted as a function of $P$, while $Z$ has been set to the equilibrium value $Z^*$. This value determines the strength of the noise coupled to prey population dynamics $dP$. Each estimation has been performed on the nonoverlapping window marked by the respective color on the bar in (a). The black dashed line represents a linear fit of all shown functions $(BB^\top)_{1,1}(P, Z^*, t)$ with respect to their dependence on $P$. All estimations were performed on running windows of length 100 time units, meaning 5000 data points at sampling rate $\Delta t = 2 \times 10^{-2}$.

The statistical quality of the estimator $\widehat{\lambda}$ with respect to its distribution width is similar to that of the conventional indicators and sufficiently good to detect CSD with a low false negative rate. This can be argued by checking that the confidence intervals at the beginning and the end of the estimator time series do not overlap (see also Appendix B for a more in-depth discussion). The deviation of the sample mean from the theoretical value of $\lambda$ [black and red curves in Fig. 1(d), respectively] originates in the nonzero time step $\Delta t = 10^{-1}$. While this does not affect the trend of the estimator time series, it decreases the signal-to-noise ratio.

### B. Predator-prey model with multiplicative noise

Following the work of Bengfort *et al.* [42], we examine a predator-prey model for oceanic plankton populations (see Appendix C for details). Since Bengfort *et al.* consider this model under the assumption of no external disturbances in the form of noise, we adopt a noise model from a related study [43]. Because environmental variability influences population growth rates, a multiplicative noise term is commonly employed [44–47]. This motivates us to investigate the following system of SDEs:

$$dP_t = \xi^{-1}\left[rP_t\left(1 - \frac{P_t}{K(\text{turb})}\right) - \frac{aP_t^2}{h(\text{turb})^2 + P_t^2}Z_t\right]dt$$
$$+ \xi^{-1/2}\sigma_P P_t dW_t^P, \quad (8)$$

$$dZ_t = \left(\frac{aP_t^2}{h(\text{turb})^2 + P_t^2}Z_t - mZ_t^2\right)dt + \sigma_Z Z_t dW_t^Z, \quad (9)$$

under the external forcing of ocean turbulence (turb). Due to the quadratic mortality term $mZ^2$ of the predator population $Z$, this system can exhibit multiple stable equilibria and bifurcations as shown in Fig. 2(a). The two white noise terms are assumed to be independent.

Here, we examine the performance of the conventional and the newly proposed CSD indicators as the system approaches the saddle-node bifurcation that annihilates the populated stable state. Figure 2(a) shows sample paths for the predator and prey populations along with the stable and unstable equilibria of the prey population $P$ as implied by the parameter turb at time $t$.

The most common approach for the assessment of CSD in a multidimensional system such as this one is to first reduce the system to one dimension [48–50]. The center manifold theorem implies that close to a critical bifurcation, the direction of lowest stability, i.e., the direction dominating the recovery time of the system, experiences further destabilization [50–52]. PCA is commonly performed to determine a linear combination of system components that exhibits the largest variance and can thus be suspected to be of the lowest stability and therefore undergo further destabilization [15,16,19,53]. However, as can be seen in the example at hand, away from the immediate proximity of the bifurcation point, the destabilizing direction need not be the direction of lowest stability. In our example, methods relying on PCA [15] or maximum autocorrelation factor (MAF) analysis [19] would identify a direction closely aligned with the predator direction $Z$ as a potential candidate for a destabilizing direction. This is problematic, as this dimension is relatively unaffected by changes in the control parameter turb and does not exhibit CSD [gray curves in Figs. 2(b) and 2(c)].

To circumvent this issue, one should therefore perform a comprehensive stability analysis on the multidimensional time series. This is achieved by examining both eigenvalues of the local equilibrium dynamics as shown in Fig. 2(d).

The $\widehat{\lambda}$ associated with the more stable eigenvalue can be seen to substantially decrease. Note that the conventional approach of focusing on the direction of the most unstable eigenvalue would miss this destabilization. An assessment of all principal components and their respective amplitudes does not alleviate this complication, as PCA could still be influenced by the dependence of the driving noise on time and on the system state. Employing the MAF method proposed in [19,54] and monitoring all autocorrelation factors would be a good alternative, as the influence of noise is eliminated through an initial normalization of the data. Fitting a discrete-time autoregression model to the time series data may also yield satisfying results [55].

Investigating again the conventional variance and AC(1)-based CSD indicators in Figs. 2(b) and 2(c), AC(1) purports destabilization along the dimension of the prey population $P$, but the trend of the observed variance seems to indicate the opposite. This contradiction can be resolved by examining the estimated diffusion matrix $BB^\top$ and its evolution along the observed time series [see Fig. 2(e)]. A clear state-dependence with respect to the prey population size $P$ can be identified in the noise amplitude $\sqrt{(BB^\top)_{1,1}}(P, Z^*, t)$ associated with d$P$. Taken together with the observation of a declining mean state in the prey population while approaching the critical transition, it can be concluded that the decrease in variance in Fig. 2(b) was due to a reduced noise amplitude and can thus be reconciled with the increase in AC(1).

## IV. DISCUSSION

The use of variance and AC(1) in combination is a well-established and by far the most common course of action to assess whether or not a system is approaching a critical transition. In general, a positive result is considered robust when both indicators show a significant positive trend [56]. We have shown here that in the presence of time- or state-dependent noise amplitudes, the variance of the system may actually decrease in the advent of a bifurcation. It may also increase without the system approaching a bifurcation (not shown). If a monitored system shows a decreasing trend in variance alongside an increasing trend in AC(1) as it approaches a bifurcation, this would typically not be considered a robust EWS, leading to a missed alarm. We have also shown that common methods in dimension reduction can lead to missed alarms, as the destabilizing system component may not be the least stable one at the initial distance to the bifurcation.

To overcome these problems, we have proposed a method based on estimating a Langevin equation from the observed dynamics. Our approach allows us to separate the effects of possible CSD contained in the drift coefficient from changes in the noise represented by the diffusion coefficient. Thereby, inconsistent indications in variance and AC(1) with respect to CSD can be reconciled. It also allows for a more holistic investigation of multidimensional systems without further mechanistic simplifications (see Appendix D for a second example to this point, which shows in particular that the proposed method works for and identifies oscillatory multidimensional systems). We showed that our approach circumvents the above-mentioned pitfalls that a conventional CSD analysis suffers from.

The confidence intervals in Fig. 1 suggest that the statistical quality of the estimator $\widehat{\lambda}$ is of the same order as that of the conventional early warning indicators. In Appendix B we investigate the data requirements of the conventional and novel indicators in a one-dimensional system in more detail. There, variance and AC(1) can be seen to be favorable in settings of low data availability, i.e., short time series and low temporal resolution. In this context, an important caveat bears mentioning once more: While for the estimator of variance, the length of the time series is the only determining factor of convergence, the estimators for the drift and diffusion coefficients also require small sampling time steps $1 \gg \lambda \Delta t > 0$ in order for their biases to be small. In general, the estimator $\widehat{\lambda}$ proposed here still contains information about CSD even in settings of large sample time steps $\Delta t$, but the signal-to-noise ratio may prohibit its employment as a CSD indicator. Areas of application where systems are potentially susceptible to tipping and where high-frequency data are available for analysis could be electricity grids [57–59], financial markets [60–62], atmospheric circulation systems such as monsoons [63,64], ecosystems and vegetation systems such as the Amazon rainforest [65–68], ocean circulation systems [3], or ice sheets [69,70]. An analogous comparison of the indicator's performance with respect to data availability in multidimensional systems is more involved. In settings where the destabilizing direction is directly aligned with an observed dimension and this dimension is correctly identified, the conventional indicators outperform our novel one, just as seen in the one-dimensional case. There are, however, also many configurations in which the conventional methods struggle to identify CSD (see Appendix B, contrasting Figs. 4 to 5). These instances should be regarded as more generic, as they occur when the deterministic eigenbasis is not precisely aligned with the directions of noise disturbances.

Our method should be understood as a more general, reliable, and circumspect indicator of CSD compared to the widely used variance and AC(1). Our approach is appropriate in settings of generally time- and state-dependent driving noise, where the combined conventional indicators may easily fail. Moreover, the ability to examine time series in their multidimensional complexity constitutes a considerable improvement in the comprehension of the system compared to one-dimensional summary statistics.

Visit the GitHub repository KramersMoyalEWS [71] to access the code generating all figures in this manuscript.

## ACKNOWLEDGMENTS

This work has received funding from the Volkswagen Stiftung, the European Union's Horizon 2020 Research and Innovation Program under Grant Agreement No. 820970 and under the Marie Sklodowska-Curie Grant Agreement No. 956170, as well as from the Federal Ministry of Education and Research under Grant No. 01LS2001A. This is ClimTip contribution #13; the ClimTip project has received funding from the European Union's Horizon Europe Research and Innovation Program under Grant Agreement No. 101137601.

## APPENDIX A: DETAILS ON THE ESTIMATOR $\widehat{\lambda}$

Here we give a detailed description of the local stability measure $\widehat{\lambda}$ in $n$-dimensional systems. For $n = 1$, $\widehat{\lambda}$ is simply given by the negative of the slope of the estimated drift coefficient $a$ around the equilibrium and is therefore real valued. For $n > 1$, the estimation procedure returns $n$ eigenvalues of the local Jacobian matrix of the estimated drift coefficient in their algebraic multiplicity. The eigenvalues may be complex, but only the negative of the real part is considered as an indicator of local stability.

Prior to any analysis, the windowed time series is linearly detrended along all $n$ dimensions. As the assessment of the conventional CSD indicators relies purely on the fluctuations around the equilibrium state, the mean of the data is disregarded. In contrast, drift and diffusion are assessed without subtraction of the mean to retain information about the corresponding state dependence. In order to obtain numerical stability, the $n$ time series are normalized to a standard deviation of 1, with no implication on the subsequent estimation of $\widehat{\lambda}$. For each window the estimation of the function $A(\mathbf{x})$ is returned as an array of values

$$(\widehat{A}(\mathbf{x}_i))_{i=1,\ldots,M^n}, \quad (A1)$$

where $M$ is the number of evenly spaced grid points in each dimension. The $M^n$ grid points in $\mathbb{R}^n$, therefore, form a grid on the hypercube spanned by the state space explored by the time series. For each of these grid points, an estimation $\widehat{A}(\mathbf{x}_i)$ is calculated using an Epanechnikov kernel,

$$K(\mathbf{x}) = \frac{3}{4h}\left(1 - \frac{||\mathbf{x}||^2}{h^2}\right), \text{ with support } ||\mathbf{x}|| < h, \quad (A2)$$

with a kernel bandwidth $h$ of $14n/M$. Since the (biased) estimator $\widehat{A}(\mathbf{x}_i)$ will converge as the number of samples $\mathbf{X}_{k\Delta t}$ close to $\mathbf{x}_i$ tends to infinity, it is clear that those grid points with many samples in their kernel proximity will see the fastest convergence. In our setting with equilibrium dynamics around one stable equilibrium $\mathbf{x}^*$, this means that the estimations for grid points closest to $\mathbf{x}^*$ converge fastest, and the quality deteriorates for outer grid points. For this reason and in order to curtail the effects of a nonlinear drift term, we opt to only carry 50% of grid points in each dimension centered around the grid point containing $\widehat{\mathbf{x}}^* = \widehat{\mathbb{E}}[\mathbf{X}]$ to the subsequent analysis. This is to say, we select a hypercube with side lengths 50% as large as the original hypercube. In the above, three free parameters have been chosen implicitly: The number of total grid points $M$ in each dimension, the percentage $m$ of grid points to carry on either side of the estimated equilibrium $\widehat{\mathbf{x}}^*$, and the kernel bandwidth. In this study we chose $M = 50$ and $m = 50\%$, meaning that for $n = 1$, we have 25 relevant grid points for further analysis. The bandwidth is chosen as a function of $M$ and $n$, as described above. However, the performance of the estimator $\widehat{\lambda}$ is not very sensitive to small changes in these parameters.

To obtain an estimation of the local Jacobian matrix around the equilibrium point $\widehat{\mathbf{x}}^*$, we perform a multivariate ordinary linear regression between $(\widehat{A}(\mathbf{x}_i))$ and $(\mathbf{x}_i)$ over $i$, including an intercept in the design matrix. The algebraic eigenvalues of the resulting matrix are computed numerically. For $n = 1$, this procedure is equivalent to finding a best linear fit $(c - \widehat{\lambda}x_i)$ to $(\widehat{a}(x_i))$ over $i$.

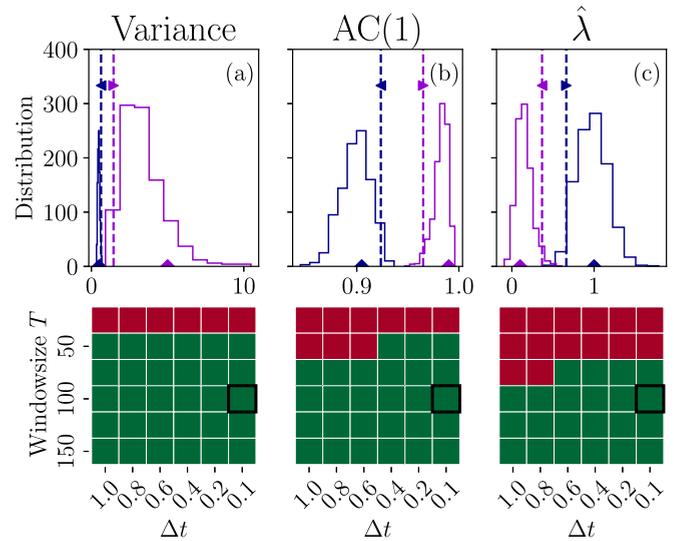

FIG. 3. Distributions of (a) variance, (b) AC(1), and (c) the estimator $\hat{\lambda}$ on 1000 samples of Ornstein–Uhlenbeck time series of length $T = 100$ each [see Eq. (B1)]. The sampling time step was set to $\Delta t = 0.1$. This estimation setting is equivalent to that chosen for the synthetic fold bifurcation example in the main text. The blue histogram shows the estimator distribution for the OU parameter $\lambda = 1$ and the violet for $\lambda = 0.1$. If the 95% confidence intervals of the two distributions do not overlap, it is highly likely that the corresponding CSD will be detected by the indicators. The separation in distribution is investigated for several choices of $T$ and $\Delta t$ in the bottom row of the figure. Green indicates no significant overlap of the estimator distributions, i.e., a successful CSD detection; red the contrary.

## APPENDIX B: ASSESSING THE STATISTICAL QUALITY OF $\hat{\lambda}$

The two applications presented in the main text demonstrate that CSD manifests itself in a substantial negative trend of the estimator $\hat{\lambda}$ when enough data are available. In this section we aim to make this statement concrete and to compare the indicator's performance to that of variance and AC(1). The assessment is based on the width of the three indicators' numerical distributions after application to one-dimensional synthetic data. The top row of Fig. 3 shows the distributions of the estimators that arise from the application of the indicators to 1000 synthetic time series generated by numerically integrating a time-homogeneous Ornstein-Uhlenbeck (OU) process:

$$dX_t = -\lambda X_t dt + dW_t. \quad (B1)$$

To analyze the performance of the estimators in a generic CSD scenario, i.e., a temporal reduction of the restoring rate, we plot their distributions for $\lambda = 1$ and $\lambda = 0.1$. If the distributions are sufficiently distinct, the indicators correctly detect this reduction of the restoring rate from 1 to 0.1 with a high likelihood.

For different choices of window lengths $T$ and time steps $\Delta t$, we check numerically whether this condition is satisfied for each estimator (bottom row of Fig. 3). At low data availability, the indicators variance and AC(1) outperform $\hat{\lambda}$. Above a window length of $T = 100$, this difference is negligible, judging by the proposed metric. The difference may also be less pronounced when performing the same test on time series generated by models with nonlinear drift or jump noise, where the state locality of our method can alleviate nonlinear effects on the far ends of the state space. Therefore, in a large range of applications, the estimator $\hat{\lambda}$ offers a statistically equally performing method of assessing CSD with the additional advantage of robustness with respect to time- and state-dependent noise, as well as settings featuring nonlinear drifts and jumps in the noise.

An analogous analysis of indicator performances on data from multidimensional systems would need to encompass many more parameter choices. A general setup for a comparison in two dimensions could be

$$d\mathbf{X}_t = O \begin{pmatrix} -\lambda_1 & 0 \\ 0 & -\lambda_2 \end{pmatrix} O^{-1} \mathbf{X}_t dt + \begin{pmatrix} \sigma_1 & c \\ 0 & \sigma_2 \end{pmatrix} d\mathbf{W}_t, \quad \text{(B2)}$$

where $O$ is an orthogonal matrix determining the alignment of the eigenspace of the inner matrix with respect to the observed processes $\mathbf{X}^{(1,2)}$ and the noise disturbances. A fold bifurcation could be generically represented in this setting by reducing $\lambda_1$ to zero while keeping $\lambda_2$, $\sigma_1$, $\sigma_2$, and $c$ fixed. The conventional indicators of CSD would be vulnerable to a time evolution in these quantities, as well as a potential time-dependent $O(t)$, while our novel method is not. In the most trivial setting of $O = \mathbf{1}$, observing $\mathbf{X}^{(1)}$ is equivalent to the one-dimensional setting above. We compare the performance of variance and AC(1) when computed on the time series of $\mathbf{X}^{(1)}$ to the most pronounced downward trend of $\hat{\lambda}$ computed on the low-dimensional time series of $\mathbf{X}$. As can be seen in Fig. 4, the performance of $\hat{\lambda}$ deteriorates in the transition from 1D to 2D data. This is likely due to the more sparse estimation of the function $A(\cdot, t)$ and the propagated uncertainties in the more unstable estimations for $DA$ and ultimately $\hat{\lambda}$. However, these enhanced complications in the estimation $\hat{\lambda}$ can in many cases be less severe than the issues arising in the indicators using variance and AC(1). If the eigenspace is rotated by 45°, i.e.,

$$O = \begin{pmatrix} 2^{-1/2} & -2^{-1/2} \\ 2^{-1/2} & 2^{-1/2} \end{pmatrix}, \quad \text{(B3)}$$

then the mixing of linear restoring rates can mask the effects of CSD in the variables $\mathbf{X}^{(1,2)}$ when considered on their own. The comprehensive analysis of eigenvalues on the two-dimensional time series, however, can still identify the underlying change in system dynamics (see Fig. 5). Since this misalignment of the eigenspace with the observed dimensions is the generic case as opposed to the specific choice of $O = \mathbf{1}$, we propose that the use of our novel indicator can be numerically favorable with respect to data requirements in a wide range of application settings.

## APPENDIX C: DETAILS ON THE PREDATOR-PREY MODEL

The specific model introduced in the main text is a modified version of the Truscott–Brindley model for ocean plankton populations originally introduced in [72]. Bengfort *et al.* [42] generalized the model by introducing the

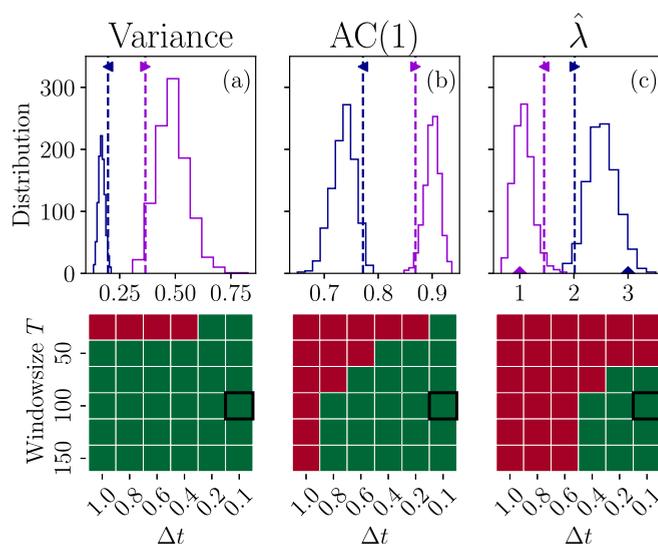

FIG. 4. Distributions of (a) variance, (b) AC(1) on 1000 samples of $\mathbf{X}^{(1)}$ time series of length $T = 100$ each computed from the 2D Ornstein–Uhlenbeck model in Eq. (B2) with $O = \mathbf{1}$. The estimator $\hat{\lambda}$ in (c) was computed on the 2D time series, and the larger of the two estimated eigenvalue real parts was investigated. The sampling time step was set to $\Delta t = 0.1$. The blue and violet histograms show the estimator distributions for the parameter $\lambda_1 = 3$ and $\lambda_1 = 1$, respectively, while $\lambda_2, \sigma_1, \sigma_2 = 1$ were identical for both histograms. If the 97.5 and 2.5 percentile of the two distributions do not overlap, it is highly likely that the corresponding CSD will be detected by the indicators. The separation in distribution is investigated for several choices of $T$ and $\Delta t$ in the bottom row of the figure. Color-coding according to Fig. 3.

environmental parameter of fluid turbulence to the system and by allowing higher powers in the mortality term of the predator population. In their nondimensionalized form the full

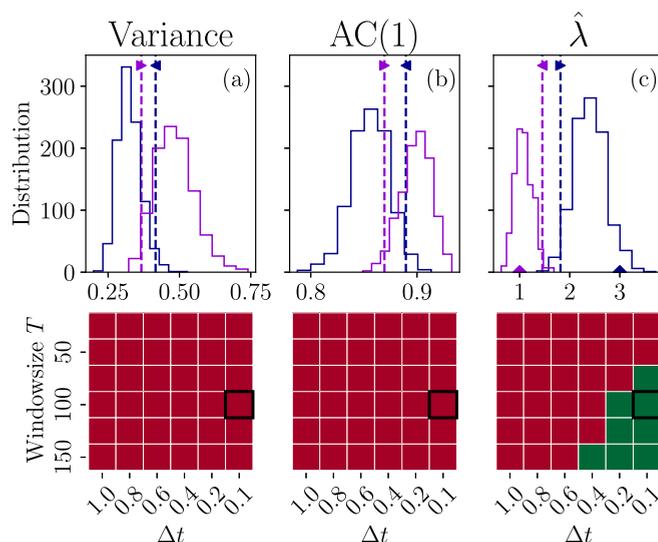

FIG. 5. Identical setup as in Fig. 4 but now for the eigenspace basis defined through the matrix $O$ given in Eq. (B3). The eigenspace is not aligned with the observed dimensions, leading to a substantial sensitivity loss in variance and AC(1).

TABLE I. Parameter values used in the simulation of plankton populations following the model in Eqs. (C5).

| Parameter | Value | Description |
|---|---|---|
| $r$ | 1 | Growth rate factor of prey $P$ |
| $a$ | 1/9 | Rate of the predator consuming the prey |
| $m$ | 0.0525 | Mortality rate of the predator |
| $\xi$ | 0.7 | Timescale separation between prey and predator evolutions |
| $h_0$ | 1/16 | Factor influencing maximal consumption at zero turbulence |
| $c_h$ | 0.88 | Linear relationship between turbulence and $h$ |
| $K_0$ | 0.7 | Carrying capacity at zero turbulence |
| $c_K$ | 0.66 | Linear relationship between turbulence and $K$ |
| $\sigma_P$ | 0.037 | Strength of noise coupling to $P$ |
| $\sigma_Z$ | 0.01 | Strength of noise coupling to $Z$ |

system equations read

$$\xi \dot{p}(t) = rp(t)\left(1 - \frac{p(t)}{K(\text{turb})}\right) - \frac{ap(t)^2}{h(\text{turb})^2 + p(t)^2}z(t), \quad \text{(C1)}$$

$$\dot{z}(t) = \frac{ap(t)^2}{h(\text{turb})^2 + p(t)^2}z(t) - mz(t)^2, \quad \text{(C2)}$$

$$K(\text{turb}) = K_0 + c_K \text{turb}, \quad \text{(C3)}$$

$$h(\text{turb}) = \frac{h_0}{1 + c_h \text{turb}}. \quad \text{(C4)}$$

The first term on the right-hand side of the prey population's evolution $\dot{p}$ is the population growth rate as determined by the relationship between the current population size and the carrying capacity $K$. The second term is the mortality rate of the prey population, which is simultaneously the growth rate of the predator population since it is assumed that all death in $p$ and growth in $z$ occurs through consumption of the former by the latter. The second term in the evolution of $z$ in the second equation is the quadratic mortality term highlighted in the main text. This term facilitates multiple stable states as opposed to the same model with a linear mortality term. The turbulence turb $\in [0, 1]$ describes the normalized strength of spatial mixing in the ocean modeled by circular eddies. All parameter values but those for $\xi$, $c_K$, and $c_h$ are adopted directly from [42] and can be found in Table I, along with a short description of their interpretation. The parameters $c_K$ and $c_h$ were increased by a factor of 2.2 each for the purposes of this study to support a bigger range of stable prey populations in the large population regime. The fundamental nature of the model remains unaltered by this change. Lastly, as described in the main text, we introduced multiplicative noise terms commonly used in the relevant literature [44–47] to model environmental impacts on the growth and mortality rates of the two populations. This leads us to the complete set of model equations:

$$dP_t = \xi^{-1}\left[rP_t\left(1 - \frac{P_t}{K(\text{turb})}\right) - \frac{aP_t^2}{h(\text{turb})^2 + P_t^2}Z_t\right]dt + \xi^{-1/2}\sigma_P P_t dW_t^P.$$

$$dZ_t = \left(\frac{aP_t^2}{h(\text{turb})^2 + P_t^2}Z_t - mZ_t^2\right)dt + \sigma_Z Z_t dW_t^Z.$$

$$K(\text{turb}) = K_0 + c_K \text{turb},$$

$$h(\text{turb}) = \frac{h_0}{1 + c_h \text{turb}}$$

$$\text{turb}(t) = 1 - \frac{7}{10}\frac{t}{T}. \quad \text{(C5)}$$

### APPENDIX D: ADDITIONAL EXAMPLE OF THE MULTIDIMENSIONAL STABILITY ANALYSIS

For fold-type bifurcations akin to the two exemplary systems in the main text, the center manifold theorem allows for an identification of a one-dimensional critical subspace close to the bifurcation point. On this subspace, the dynamics are topologically equivalent to the normal form of the fold bifurcation. However, the critical dynamics of many systems cannot be reduced to a one-dimensional subspace in such a way. This is especially relevant for systems exhibiting pronounced oscillations. Using the method outlined above, we therefore additionally assess the local stability of a system undergoing a subcritical Hopf bifurcation in normal form:

$$d\mathbf{X}_t = \begin{pmatrix} -(\mu(t) - (\mathbf{X}^{(1)})^2 - (\mathbf{X}^{(2)})^2)\mathbf{X}^{(1)} - \omega\mathbf{X}^{(2)} \\ -(\mu(t) - (\mathbf{X}^{(1)})^2 - (\mathbf{X}^{(2)})^2)\mathbf{X}^{(2)} + \omega\mathbf{X}^{(1)} \end{pmatrix} dt \quad \text{(D1)}$$

$$+ \varepsilon \begin{pmatrix} 1 & 0 \\ 0 & 1 \end{pmatrix} d\mathbf{W}_t, \quad \text{(D2)}$$

where $\omega = 1$, $\varepsilon = 0.01$, and $\mu(t)$ decreases linearly from 2 to 0.1 over the integration time of $T = 1000$. For $\mu > 0$, the origin is a stable fixed point with eigenvalues $-\mu \pm i\omega$. Furthermore, there is an unstable limit cycle with radius $\sqrt{\mu}$, and perturbations from the origin decay in the form of spirals. At $\mu = 0$, the radius of the unstable limit cycle reaches zero, and the origin turns into an unstable fixed point. The data was sampled at time steps $\Delta t = 0.1$ and analyzed in windows of length $T = 100$. The results are presented in Fig. 6. The real parts of both eigenvalues are known to be $-\mu(t)$, and the estimations track this value relatively closely. A destabilization of the equilibrium can clearly be made out. An additional insight gained via the CSD assessment through the Langevin equation approach proposed here is that the local system exhibits oscillatory dynamics, as identified by the complex eigenvalues of the Jacobian matrix.

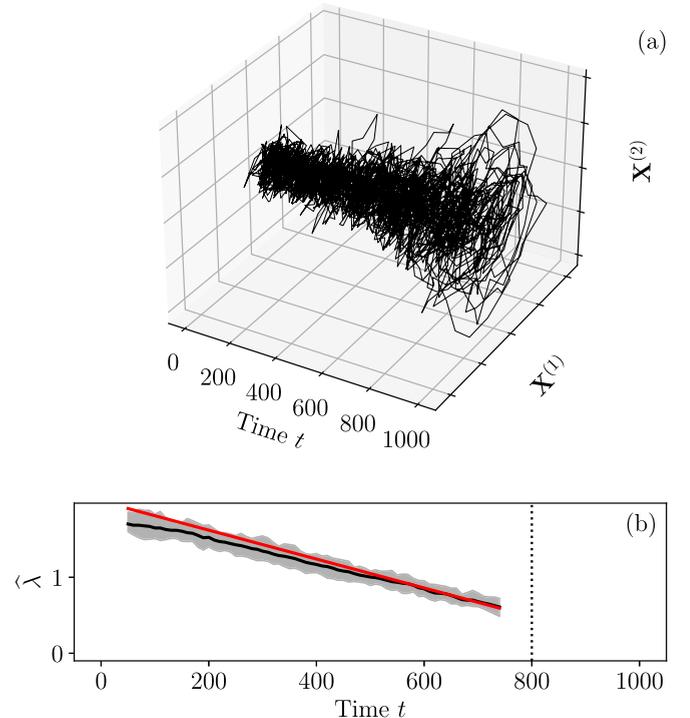

FIG. 6. Employing the new method for multidimensional stability analysis. (a) Time series data generated by a model approaching a subcritical Hopf bifurcation in normal form. (b) Estimates of the stability indicators $\lambda_k$. The thick line represents the mean of the 200 sample estimations. Since the Jacobian matrices at each point in time have complex eigenvalues corresponding to the oscillatory dynamics, the eigenvalues are conjugates and coincide in their real part. The theoretical value for the negative of this real part is plotted in red.